\documentclass[12pt]{article}

\usepackage{amsmath, epsfig, cite}
\usepackage{amssymb}
\usepackage{amsfonts}
\usepackage{latexsym}
\usepackage{float}

\newtheorem{thm}{Theorem}[section]

\numberwithin{equation}{section}

\newcommand{\qed}{{\hfill$\square$}\medskip}
\UseRawInputEncoding
\begin{document}


\begin{center}
{\Large\bf Combinatorial proofs of several partition identities of\\[7pt]
Andrews and Merca}
\end{center}

\vskip 2mm \centerline{Ji-Cai Liu and Huan Liu}
\begin{center}
{\footnotesize Department of Mathematics, Wenzhou University, Wenzhou 325035, PR China\\
{\tt jcliu2016@gmail.com, 23450829016@stu.wzu.edu.cn} }
\end{center}


\vskip 0.7cm \noindent{\bf Abstract.}
Based on a bijection due to Fu and Tang, we provide combinatorial proofs of several partition identities of Andrews and Merca. We also introduce two weights for partitions to extend one of these identities.


\vskip 2mm
\noindent{\it MR Subject Classifications}: 11P84, 05A17, 05A19
\section{Introduction}
A partition of a positive integer $n$ is a finite nonincreasing sequence of positive integers $\lambda_1,\lambda_2,\cdots,\lambda_r$ such that $\sum_{i=1}^r\lambda_i=n$. A partition of $n$ can be rewritten as
\begin{align*}
n=t_1+2t_2+\cdots+nt_n,
\end{align*}
where each positive integer $i$ appears $t_i$ times in the partition.

Andrews and Merca \cite[Theorems 1.1 and 1.4]{am-jcta-2024} established the following two interesting results:
\begin{align}
&\sum_{t_1+2t_2+\cdots+nt_n=n}(\pm 1)^{t_2+t_4+t_6+\cdots}\left(t_2+2t_4+3t_6+\cdots\right)\notag\\
&=\sum_{t_1+2t_2+\cdots+nt_n=n}(\pm 1)^{\lfloor t_1/2\rfloor+\lfloor t_2/2\rfloor+\cdots +\lfloor t_n/2\rfloor}\left(\left\lfloor t_1/2 \right\rfloor+2\left\lfloor t_2/2 \right\rfloor+\cdots+n\left\lfloor t_n/2 \right\rfloor\right),\label{am-1}
\end{align}
and
\begin{align}
&\sum_{t_1+2t_2+\cdots+nt_n=n}(\pm 1)^{t_2+t_4+t_6+\cdots}\left(t_1+3t_3+5t_5+\cdots\right)\notag\\
&=\sum_{t_1+2t_2+\cdots+nt_n=n}(\pm 1)^{\lfloor t_1/2\rfloor+\lfloor t_2/2\rfloor+\cdots +\lfloor t_n/2\rfloor}\left(\widehat{t_1}+2\widehat{t_2}+\cdots+n \widehat{t_n}\right),\label{am-2}
\end{align}
where $\lfloor x\rfloor$ denotes the integral part of real $x$ and $\widehat{t}$ denotes the residue of $t$ modulo $2$.

The proof of the results \eqref{am-1} and \eqref{am-2} relies on generating functions. ``Combinatorial proofs of these results would be very interesting" as described by Andrews and Merca \cite{am-jcta-2024}.

Recently, Merca \cite[Theorems 1.3 and 1.7]{merca-racsam-2024} established more integer partition identities as follows.
\begin{thm}\label{t-1}
Let $z$ be a complex number.
For all positive integers $m$ and $n$, we have
\begin{align}
&\sum_{t_1+2t_2+\cdots+nt_n=n}(t_m\pm 2^zt_{2m}+3^zt_{3m}\pm\cdots)\notag\\[5pt]
&=\sum_{t_1+2t_2+\cdots+nt_n=n}\left(\left\lfloor \frac{t_1}{m}\right\rfloor\pm 2^z\left\lfloor \frac{t_2}{m}\right\rfloor+3^z\left\lfloor \frac{t_3}{m}\right\rfloor+\cdots\right),\label{m-1}\\[10pt]
&\sum_{t_1+2t_2+\cdots+nt_n=n}(t_m\pm 2^zt_{2m}+3^zt_{3m}\pm\cdots)\notag\\[5pt]
&=\sum_{t_1+2t_2+\cdots+nt_n=n}\left(\left\lfloor \frac{t_1}{m}\right\rfloor\pm 2^z\left\lfloor \frac{t_1}{2m}\right\rfloor+3^z\left\lfloor \frac{t_1}{3m}\right\rfloor+\cdots\right),\label{m-2}
\end{align}
and
\begin{align}
&\sum_{t_1+2t_2+\cdots+nt_n=n}(-1)^{t_m+t_{2m}+t_{3m}+\cdots}(t_m\pm 2^zt_{2m}+3^zt_{3m}\pm\cdots)\notag\\[5pt]
&=\sum_{t_1+2t_2+\cdots+nt_n=n}(-1)^{\lfloor t_1/m\rfloor+\lfloor t_2/m\rfloor+\lfloor t_3/m\rfloor+\cdots}\left(\left\lfloor \frac{t_1}{m}\right\rfloor\pm 2^z\left\lfloor \frac{t_2}{m}\right\rfloor+3^z\left\lfloor \frac{t_3}{m}\right\rfloor+\cdots\right).\label{m-3}
\end{align}
\end{thm}

The aim of the paper is to give combinatorial proofs of \eqref{am-1}--\eqref{m-3}. In fact, we prove \eqref{am-1} and \eqref{am-2} in the following general forms.
\begin{thm}\label{t-2}
For all positive integers $m$ and $n$, we have
\begin{align}
&\sum_{t_1+2t_2+\cdots+nt_n=n}(\pm 1)^{t_m+t_{2m}+t_{3m}+\cdots}\left(t_m+2t_{2m}+3t_{3m}+\cdots\right)\notag\\
&=\sum_{t_1+2t_2+\cdots+nt_n=n}(\pm 1)^{\lfloor t_1/m\rfloor+\lfloor t_2/m\rfloor+\cdots +\lfloor t_n/m\rfloor}\left(\left\lfloor t_1/m \right\rfloor+2\left\lfloor t_2/m \right\rfloor+\cdots+n\left\lfloor t_n/m \right\rfloor\right),\label{a-new-1}
\end{align}
and
\begin{align}
&\sum_{t_1+2t_2+\cdots+nt_n=n}(\pm 1)^{t_m+t_{2m}+t_{3m}+\cdots}\sum_{i\ge 0}\sum_{j=1}^{m-1}(im+j)t_{im+j}\notag\\
&=\sum_{t_1+2t_2+\cdots+nt_n=n}(\pm 1)^{\lfloor t_1/m\rfloor+\lfloor t_2/m\rfloor+\cdots +\lfloor t_n/m\rfloor}\left(\langle t_1\rangle_m+2\langle t_2\rangle_m+\cdots+n \langle t_n\rangle_m\right),
\label{a-new-2}
\end{align}
where $\langle x\rangle_m$ denotes the least non-negative integer $r$ with $x\equiv r \pmod{m}$.
\end{thm}

The rest of the paper is organized as follows. In the next section, we first recall a bijection due to Fu and Tang \cite{ft-ms-2017}, from which we establish an equidistributed result on the set of all partitions of $n$. By using this equidistributed result, we prove Theorems \ref{t-1} and \ref{t-2} in Sections 3 and 4, respectively. In the last section, we extend \eqref{m-2} by introducing two weights for partitions.

\section{A bijection of Fu and Tang}
For an integer $m\ge 2$, let $O_m(n)$ denote the set of partitions of $n$ in which no part is a multiple of $m$, and $D_m(n)$ denote the set of partitions of $n$ in which all parts appear less than $m$ times.
There is a beautiful bijection between $O_m(n)$ and $D_m(n)$ due to Glaisher \cite{lehmer-bams-1946}.

{\noindent\it The map $\varphi^{od}_{m}$ from $O_m(n)$ to $D_m(n)$:} If the partition contains $m$ copies of the same part, then we merge the $m$ parts into one part of $m$-fold size. We repeat this procedure until all parts appear less than $m$ times.

{\noindent\it The map $\varphi^{do}_{m}$ from $D_m(n)$ to $O_m(n)$:} If the partition contains a part which is a multiple of $m$, then we split this part into $m$ equal parts. We repeat this procedure until no part is a multiple of $m$.

For an integer $m\ge 2$, let $E_m(n)$ denote the set of partitions of $n$ in which all parts are multiples of $m$, and $N_m(n)$ denote the set of partitions of $n$ in which the multiplicity of each part is a multiple of $m$. There also exists a bijection between $E_m(n)$ and $N_m(n)$.

{\noindent\it The map $\psi^{en}_{m}$ from $E_m(n)$ to $N_m(n)$:}
We split a part which is a multiple of $m$ into $m$ equal parts.

{\noindent\it The map $\psi^{ne}_{m}$ from $N_m(n)$ to $E_m(n)$:}
We merge $m$ equal parts into one part of $m$-fold size.

A partition of a positive integer $n$: $\lambda_1,\lambda_2,\cdots,\lambda_r$ can be represented as the multiset $M=\{\lambda_1,\lambda_2,\cdots,\lambda_r\}$. The multiset $M$ possesses a unique decomposition into two multisets:
\begin{align*}
M=M^{o}_{m}\cup M^{e}_{m}
\end{align*}
such that no element in $M^{o}_{m}$ is a multiple of $m$ and all elements in $M^{e}_{m}$ are multiples of $m$. For example,
\begin{align*}
\{1,1,1,2,2,2,3,3,4,4,5,6\}=\{1,1,1,2,2,2,4,4,5\}\cup \{3,3,6\},
\end{align*}
for $m=3$.

The multiset $M$ also possesses a unique decomposition into two multisets:
\begin{align*}
M=M^{d}_{m}\cup M^{n}_{m}
\end{align*}
such that all elements in $M^{d}_{m}$ appear less than $m$ times and the multiplicity of each element in $M^{n}_{m}$ is a multiple of $m$. For example,
\begin{align*}
\{1,1,1,1,1,1,2,2,2,3,4,4,5,6\}=\{3,4,4,5,6\}\cup \{1,1,1,1,1,1,2,2,2\},
\end{align*}
for $m=3$.

We recall a bijection $\sigma_m:\mathcal{P}_n\to \mathcal{P}_n$ due to Fu and Tang \cite{ft-ms-2017}.
Let $n$ and $m$ be positive integers with $m\ge 2$.
For any $M\in \mathcal{P}_n$, the map $\sigma_m$ is defined by
\begin{align*}
\sigma_{m}\left(M\right)=\varphi^{od}_{m}\left(M^{o}_{m}\right)\cup \psi^{en}_{m}\left(M^{e}_{m}\right).
\end{align*}
For example,
\begin{align*}
\sigma_{3}:\quad&\{1,1,1,2,2,2,3,3,4,4,5,6\}=\{1,1,1,2,2,2,4,4,5\}\cup \{3,3,6\}\\
&\to \{3,4,4,5,6\}\cup \{1,1,1,1,1,1,2,2,2\}=\{1,1,1,1,1,1,2,2,2,3,4,4,5,6\}.
\end{align*}
Note that all elements in $\varphi^{od}_{m}(M^{o}_{m})$ appear less than $m$ times
and the multiplicity of each element in $\psi^{en}_{m}(M^{e}_{m})$ is a multiple of $m$.

The inverse of the map $\sigma_{m}$ is given by
\begin{align*}
\sigma_{m}^{-1}\left(M\right)=\varphi^{do}_{m}\left(M^{d}_{m}\right)\cup \psi^{ne}_{m}\left(M^{n}_{m}\right).
\end{align*}
For example,
\begin{align*}
\sigma_{3}^{-1}:\quad &\{1,1,1,1,1,1,2,2,2,3,4,4,5,6\}=\{3,4,4,5,6\}\cup \{1,1,1,1,1,1,2,2,2\}\\
&\to \{1,1,1,2,2,2,4,4,5\}\cup\{3,3,6\}=\{1,1,1,2,2,2,3,3,4,4,5,6\}.
\end{align*}

Let $m,n$ and $k$ be positive integers.
For a partition $\lambda=M^{o}_{m}\cup M^{d}_{m}$: $n=t_1+2t_2+\cdots+nt_n$, we introduce the following four partition statistics:
\begin{align*}
\alpha_{k,m}(\lambda)&=t_{km},\\[5pt]
\alpha'_{k,m}(\lambda)&=\lfloor t_k/m\rfloor,\\[5pt]
\gamma_m(\lambda)&=\text{the sum of the elements in $M^{o}_{m}$}\\[5pt]
&=t_1+\cdots+(m-1)t_{m-1}+(m+1)t_{m+1}+\cdots+(2m-1)t_{2m-1}+\cdots,\\[5pt]
\gamma'_m(\lambda)&=\text{the sum of the elements in $M^{d}_{m}$}\\[5pt]
&=\langle t_1\rangle_m+2\langle t_2\rangle_m+\cdots+n \langle t_n\rangle_m.
\end{align*}

Since $\sigma_{m}$ maps $M^{o}_{m}$ onto $M^{d}_{m}$, we have
$\gamma_m(\lambda)=\gamma'_m\left(\sigma_m(\lambda)\right)$ for all partitions $\lambda\in \mathcal{P}_n$. For $km\in M^{e}_{m}$ with $1\le k\le n$, the map $\sigma_{m}$ splits the part $km$ into $m$ equal parts $k$ and $\sigma_{m}$ maps $M^{o}_{m}$ onto the multiset in which all elements appear less than $m$ times, and so $\alpha_{k,m}(\lambda)=\alpha'_{k,m}\left(\sigma_m(\lambda)\right)$ for all partitions $\lambda\in \mathcal{P}_n$. Based on these two facts, we establish the following result.

\begin{thm}
For all positive integers $m$ and $n$, we have
\begin{align}
\sum_{\lambda\in \mathcal{P}_n} x_1^{\alpha_{1,m}(\lambda)}x_2^{\alpha_{2,m}(\lambda)} \cdots x_n^{\alpha_{n,m}(\lambda)} z^{\gamma_m(\lambda)}=\sum_{\lambda\in \mathcal{P}_n}
x_1^{\alpha'_{1,m}(\lambda)}x_2^{\alpha'_{2,m}(\lambda)} \cdots x_n^{\alpha'_{n,m}(\lambda)}
z^{\gamma'_m(\lambda)},\label{b-new-1}
\end{align}
where $\mathcal{P}_n$ denotes the set of partitions of $n$.
\end{thm}

Note that $\alpha_{k,m}(\lambda)=\alpha'_{k,m}(\lambda)$ and $\gamma_m(\lambda)=\gamma'_m(\lambda)$
for $m=1$. Then the case $m=1$ of \eqref{b-new-1} clearly holds.

\section{Proof of Theorem \ref{t-1}}
By \eqref{b-new-1}, we have
\begin{align*}
&x_k\frac{\text{d}}{\text{d}x_k}\left( \sum_{\lambda\in \mathcal{P}_n} x_1^{\alpha_{1,m}(\lambda)}x_2^{\alpha_{2,m}(\lambda)} \cdots x_n^{\alpha_{n,m}(\lambda)} \right){\Bigg |}_{x_1=x_2=\cdots=x_n=x}\\[10pt]
&=x_k\frac{\text{d}}{\text{d}x_k}\left(\sum_{\lambda\in \mathcal{P}_n}
x_1^{\alpha'_{1,m}(\lambda)}x_2^{\alpha'_{2,m}(\lambda)} \cdots x_n^{\alpha'_{n,m}(\lambda)} \right){\Bigg |}_{x_1=x_2=\cdots=x_n=x},
\end{align*}
which implies that
\begin{align}
\sum_{\lambda\in \mathcal{P}_n} \alpha_{k,m}(\lambda) x^{\alpha_{m}(\lambda)}
=\sum_{\lambda\in \mathcal{P}_n} \alpha'_{k,m}(\lambda) x^{\alpha'_{m}(\lambda)},\label{c-new-1}
\end{align}
where $\alpha_m(\lambda)=\sum_{k\ge 1}\alpha_{k,m}(\lambda)$ and $\alpha'_m(\lambda)=\sum_{k\ge 1}\alpha'_{k,m}(\lambda)$.

Letting $x=1$ in \eqref{c-new-1} gives
\begin{align}
\sum_{\lambda\in \mathcal{P}_n} t_{km}=\sum_{\lambda\in \mathcal{P}_n} \lfloor t_k/m\rfloor.\label{c-new-2}
\end{align}
Letting $k\to 1$ and $m\to km$ in \eqref{c-new-2}, we obtain
\begin{align}
\sum_{\lambda\in \mathcal{P}_n} t_{km}=\sum_{\lambda\in \mathcal{P}_n} \lfloor t_1/km\rfloor.
\label{c-new-3}
\end{align}
Then the proof of \eqref{m-1} and \eqref{m-2} follows from \eqref{c-new-2} and \eqref{c-new-3}.

Letting $x=-1$ in \eqref{c-new-1} gives
\begin{align}
\sum_{\lambda\in \mathcal{P}_n}(-1)^{t_m+t_{2m}+t_{3m}+\cdots}t_{km}
=\sum_{\lambda\in \mathcal{P}_n} (-1)^{\lfloor t_1/m\rfloor+\lfloor t_2/m\rfloor+\lfloor t_3/m\rfloor+\cdots} \lfloor t_k/m\rfloor.\label{c-new-4}
\end{align}
Then the proof of \eqref{m-3} follows from \eqref{c-new-4}.

\section{Proof of Theorem \ref{t-2}}
Letting $x_i\to xy^i$ for $i=1,2,\cdots,n$ in \eqref{b-new-1} gives
\begin{align*}
&\sum_{\lambda\in \mathcal{P}_n} x^{\alpha_{1,m}(\lambda)+\alpha_{2,m}(\lambda)+\cdots+\alpha_{n,m}(\lambda)}
y^{\alpha_{1,m}(\lambda)+2\alpha_{2,m}(\lambda)+\cdots+n\alpha_{n,m}(\lambda)}z^{\gamma_m(\lambda)}\\[5pt]
&=\sum_{\lambda\in \mathcal{P}_n} x^{\alpha'_{1,m}(\lambda)+\alpha'_{2,m}(\lambda)+\cdots+\alpha'_{n,m}(\lambda)}
y^{\alpha'_{1,m}(\lambda)+2\alpha'_{2,m}(\lambda)+\cdots+n\alpha'_{n,m}(\lambda)}z^{\gamma'_m(\lambda)},
\end{align*}
which is
\begin{align}
\sum_{\lambda\in \mathcal{P}_n} x^{\alpha_m(\lambda)}y^{\beta_m(\lambda)}z^{\gamma_m(\lambda)}=\sum_{\lambda\in \mathcal{P}_n} x^{\alpha'_m(\lambda)}y^{\beta'_m(\lambda)}z^{\gamma'_m(\lambda)},\label{d-new-1}
\end{align}
where $\beta_m(\lambda)=\sum_{k\ge 1}k\alpha_{k,m}(\lambda)$ and
$\beta'_m(\lambda)=\sum_{k\ge 1}k\alpha'_{k,m}(\lambda)$.

By \eqref{d-new-1}, we have
\begin{align}
\frac{\text{d}}{\text{d}y}\left(\sum_{\lambda\in \mathcal{P}_n} x^{\alpha_m(\lambda)}y^{\beta_m(\lambda)}z^{\gamma_m(\lambda)}\right){\Bigg |}_{\substack{x=\pm 1\\y=z=1}}=\frac{\text{d}}{\text{d}y}\left(\sum_{\lambda\in \mathcal{P}_n} x^{\alpha'_m(\lambda)}y^{\beta'_m(\lambda)}z^{\gamma'_m(\lambda)}\right){\Bigg |}_{\substack{x=\pm 1\\y=z=1}},\label{d-new-2}
\end{align}
and
\begin{align}
\frac{\text{d}}{\text{d}z}\left(\sum_{\lambda\in \mathcal{P}_n} x^{\alpha_m(\lambda)}y^{\beta_m(\lambda)}z^{\gamma_m(\lambda)}\right){\Bigg |}_{\substack{x=\pm 1\\y=z=1}}=\frac{\text{d}}{\text{d}z}\left(\sum_{\lambda\in \mathcal{P}_n} x^{\alpha'_m(\lambda)}y^{\beta'_m(\lambda)}z^{\gamma'_m(\lambda)}\right){\Bigg |}_{\substack{x=\pm 1\\y=z=1}}.\label{d-new-3}
\end{align}
Then the proof of \eqref{a-new-1} and \eqref{a-new-2} follows from \eqref{d-new-2} and \eqref{d-new-3}.

\section{An extension of \eqref{m-2}}
Let $\widetilde{\mathcal{P}}(n)$ denote the set of ways of writing the integer $n$ as a sum of positive integers in non-increasing order in which at most one integer may be drawn a wavy line (we always draw the first occurrence of an integer). For example,
\begin{align*}
&\widetilde{\mathcal{P}}(5)\\[5pt]
&=\left\{(1,1,1,1),(2,1,1),(2,2),(3,1),(4),(\widetilde{1},1,1,1),
(\widetilde{2},1,1),(2,\widetilde{1},1),(\widetilde{2},2),(\widetilde{3},1),(3,\widetilde{1}),(\widetilde{4})\right\}.
\end{align*}

We define two weights for $\lambda\in \widetilde{\mathcal{P}}(n)$ as follows:
\begin{align*}
W(\lambda)=x_{\lambda_1}+x_{\lambda_2}+\cdots +x_{\lambda_s},
\end{align*}
where $\lambda_1,\cdots, \lambda_s$ are all the non-drawn parts in $\lambda$, and
\begin{align*}
\widetilde{W}(\lambda)=\left\lfloor\frac{\lambda_i}{1}\right\rfloor x_1+\left\lfloor\frac{\lambda_i}{2}\right\rfloor x_2+\left\lfloor\frac{\lambda_i}{3}\right\rfloor x_3+\cdots,
\end{align*}
where $\lambda_i$ is the drawn part in $\lambda$. For example,
\begin{align*}
&W\left(2,\widetilde{1},1\right)=x_1+x_2,\\[5pt]
&W\left(\widetilde{4}\right)=0,\\[5pt]
&\widetilde{W}\left(\widetilde{3},1\right)=3x_1+x_2+x_3,\\[5pt]
&\widetilde{W}\left(3,1\right)=0.
\end{align*}

We extend \eqref{m-2} by establishing the following result.
\begin{thm}
For all positive integers $n$, we have
\begin{align}
\sum_{\lambda\in \widetilde{\mathcal{P}}(n)}W(\lambda)=\sum_{\lambda\in \widetilde{\mathcal{P}}(n)}\widetilde{W}(\lambda).\label{e-new-3}
\end{align}
\end{thm}

{\noindent\it Proof.}
Note that the coefficient of $x_j$ on the left-hand side of \eqref{e-new-3} equals the number of times that the non-drawn part $j$ appears in $\widetilde{\mathcal{P}}(n)$.
Since the number of partitions of $i$ in which the part $j$ appears at least $l$ times equals $p(i-lj)$, we have the number of times that $j$ appears in all partitions of $i$ equals
\begin{align*}
\sum_{l=1}^{\lfloor i/j\rfloor} p(i-lj).
\end{align*}
Then the number of times that the non-drawn part $j$ appears in the subset of
$\widetilde{\mathcal{P}}(n)$ in which $k$ is drawn a wavy line equals
\begin{align*}
\sum_{l=1}^{\lfloor (n-k)/j\rfloor} p(n-k-lj).
\end{align*}
It follows that the number of times that the non-drawn part $j$ appears in $\widetilde{\mathcal{P}}(n)$ equals
\begin{align*}
&\sum_{k=0}^n\sum_{l=1}^{\lfloor (n-k)/j\rfloor} p(n-k-lj)\\[5pt]
&=\sum_{i=0}^n \sum_{l=1}^{\lfloor i/j\rfloor } p(i-lj)\\[5pt]
&=\sum_{l=1}^{ \lfloor n/j\rfloor }\sum_{i=lj}^n  p(i-lj)\\[5pt]
&=\sum_{k=0}^{n-j}\left\lfloor \frac{n-k}{j}\right\rfloor p(k)\\[5pt]
&=\sum_{k=j}^{n} \left\lfloor \frac{k}{j}\right\rfloor p(n-k).
\end{align*}
We conclude that the coefficient of $x_j$ on the left-hand side of \eqref{e-new-3} equals
\begin{align*}
\sum_{k=j}^{n} \left\lfloor \frac{k}{j}\right\rfloor p(n-k).
\end{align*}

Note that the number of times that the drawn part $k$ appears in $\widetilde{\mathcal{P}}(n)$ equals
$p(n-k)$. On the right-hand side of \eqref{e-new-3}, all the drawn parts $k$'s in $\widetilde{\mathcal{P}}(n)$ produce
\begin{align*}
\sum_{j\ge 1}p(n-k)\left\lfloor \frac{k}{j}\right\rfloor x_j.
\end{align*}
It follows that the coefficient of $x_j$ on the right-hand side of \eqref{e-new-3} equals
\begin{align*}
\sum_{k=j}^{n} \left\lfloor \frac{k}{j}\right\rfloor p(n-k).
\end{align*}
Then the coefficient of $x_j$ on both sides of \eqref{e-new-3} are equal for $j\ge 1$.
This completes the proof of \eqref{e-new-3}.
\qed

To prove \eqref{m-2}, it suffices to show that
\begin{align}
&\sum_{t_1+2t_2+\cdots+nt_n=n}(x_1t_m+x_2t_{2m}+x_3t_{3m}+\cdots)\notag\\[5pt]
&=\sum_{t_1+2t_2+\cdots+nt_n=n}\left(x_1\left\lfloor \frac{t_1}{m}\right\rfloor+x_2\left\lfloor \frac{t_1}{2m}\right\rfloor+x_3\left\lfloor \frac{t_1}{3m}\right\rfloor+\cdots\right).\label{e-new-1}
\end{align}

Let $p_m(n)$ be the number of the partitions of $n$ in which no part is divisible by $m$
and the multiplicity of $1$ is less than $m$. For the sake of convenience, let $p_m(0)=1$.
Note that
\begin{align}
&\sum_{t_1+2t_2+\cdots+nt_n=n}(x_1t_m+x_2t_{2m}+x_3t_{3m}+\cdots)\notag\\[5pt]
&=\sum_{l=0}^{n}p_m(n-l)\sum_{\substack{t_1+mt_{m}+2mt_{2m}+3mt_{3m}+\cdots=l\\[5pt]
m|t_1}}(x_1t_m+x_2t_{2m}+x_3t_{3m}+\cdots),\label{e-new-4}
\end{align}
and
\begin{align}
&\sum_{t_1+2t_2+\cdots+nt_n=n}\left(x_1\left\lfloor \frac{t_1}{m}\right\rfloor+x_2\left\lfloor \frac{t_1}{2m}\right\rfloor+x_3\left\lfloor \frac{t_1}{3m}\right\rfloor+\cdots\right)\notag\\[5pt]
&=\sum_{l=0}^{n}p_m(n-l)\sum_{\substack{t_1+mt_{m}+2mt_{2m}+3mt_{3m}+\cdots=l\\[5pt] m|t_1}}\left(x_1\left\lfloor \frac{t_1}{m}\right\rfloor+x_2\left\lfloor \frac{t_1}{2m}\right\rfloor+x_3\left\lfloor \frac{t_1}{3m}\right\rfloor+\cdots\right).\label{e-new-5}
\end{align}

Combining \eqref{e-new-4} and \eqref{e-new-5}, to prove \eqref{e-new-1}, it suffices to show that
for any nonnegative integer $l$,
\begin{align}
&\sum_{\substack{t_1+mt_{m}+2mt_{2m}+3mt_{3m}+\cdots=l\\[5pt]
m|t_1}}(x_1t_m+x_2t_{2m}+x_3t_{3m}+\cdots)\notag\\[5pt]
&=\sum_{\substack{t_1+mt_{m}+2mt_{2m}+3mt_{3m}+\cdots=l\\[5pt] m|t_1}}\left(x_1\left\lfloor \frac{t_1}{m}\right\rfloor+x_2\left\lfloor \frac{t_1}{2m}\right\rfloor+x_3\left\lfloor \frac{t_1}{3m}\right\rfloor+\cdots\right),\label{e-new-2}
\end{align}
which is stronger than \eqref{e-new-1}.
Note that \eqref{e-new-2} is equivalent to \eqref{e-new-3}.

\vskip 5mm \noindent{\bf Acknowledgments.}
The first author was supported by the National Natural Science Foundation of China (grant 12171370).

\end{document}